\documentclass{article}
\usepackage{latexsym, amssymb, amsmath, amsthm, graphics}
\usepackage{url}
\usepackage{eucal}

\title{ \bf The Jacobian Conjecture as a problem in combinatorics} \author{David Wright}

\begin{document}
    \maketitle
\begin{abstract}
    The Jacobian Conjecture has been reduced to the symmetric
    homogeneous case.  In this paper we give an inversion formula for
    the symmetric case and relate it to a combinatoric structure
    called the Grossman-Larson Algebra.  We use these tools to prove
    the symmetric Jacobian Conjecture for the case $F=X-H$ with $H$
    homogeneous and $JH^{3}=0$.  Other special results are also
    derived.  We pose a combinatorial statement which would give
    a complete proof the Jacobian Conjecture.
    \end{abstract}

    \theoremstyle{definition} 
    \newtheorem{rema}{Remark}[section]
    \newtheorem{questions}[rema]{Questions}
    \newtheorem{assertion}[rema]{Assertion}
    \theoremstyle{plain} 
    \newtheorem{propo}[rema]{Proposition}
    \newtheorem{theo}[rema]{Theorem} 
    \newtheorem{conj}[rema]{Conjecture}
    \newtheorem{quest}[rema]{Question} 
    \theoremstyle{definition}
    \newtheorem{defi}[rema]{Definition}
    \theoremstyle{plain}
    \newtheorem{lemma}[rema]{Lemma} 
    \newtheorem{corol}[rema]{Corollary}
    \newtheorem{exam}[rema]{Example} 
    \newtheorem{rmk}[rema]{Remark}
    \newcommand{\del}{\triangledown}
    \newcommand{\nno}{\nonumber} \newcommand{\lbar}{\big\vert}
    \newcommand{\mbar}{\mbox{\large $\vert$}} \newcommand{\p}{\partial}
    \newcommand{\dps}{\displaystyle} \newcommand{\bra}{\langle}
    \newcommand{\ket}{\rangle} \newcommand{\kr}{\mbox{\rm Ker}\ }
    \newcommand{\res}{\mbox{\rm Res}} \renewcommand{\hom}{\mbox{\rm Hom}}
    \newcommand{\pf}{{\it Proof:}\hspace{2ex}}
    \newcommand{\epf}{\hspace{2em}$\Box$}
    \newcommand{\epfv}{\hspace{1em}$\Box$\vspace{1em}}
    \newcommand{\nord}{\mbox{\scriptsize ${\circ\atop\circ}$}}
    \newcommand{\wt}{\mbox{\rm wt}\ } \newcommand{\clr}{\mbox{\rm clr}\ }
    \newcommand{\ideg}{\mbox{\rm Ideg}\ } \newcommand{\GC}{G_{\mathbb C}}
    \newcommand{\gC}{{\mathfrak g}_{\mathbb C}}
    \newcommand{\hatC}{\widehat {\mathbb C}} \newcommand{\bC}{{\mathbb C}}
    \newcommand{\bZ}{{\mathbb Z}} \newcommand{\bQ}{{\mathbb Q}}
    \newcommand{\bR}{{\mathbb R}} \newcommand{\bN}{{\mathbb N}}
    \newcommand{\bT}{{\mathbb T}} \newcommand{\fg}{{\mathfrak g}}
    \newcommand{\fgC}{{\mathfrak g}_{\bC}} \newcommand{\Em}{{\cal E}(m)}
    \newcommand{\cD}{\mathcal D} \newcommand{\cP}{\mathcal P}
    \newcommand{\cC}{\mathcal C} \newcommand{\cS}{\mathcal S}
    \newcommand{\EGC}{{\cal E}(\GC)}
    \newcommand{\cLGC}{\widetilde{L}_{an}\GC}
    \newcommand{\LGC}{{L}_{an}\GC} \newcommand{\BQ}{\begin{eqnarray}}
    \newcommand{\EQ}{\end{eqnarray}} \newcommand{\BQn}{\begin{eqnarray*}}
    \newcommand{\EQn}{\end{eqnarray*}} \newcommand{\wtilde}{\widetilde}
    \newcommand{\Hol}{\mbox{Hol}} \newcommand{\Hom}{\mbox{Hom}}
    \newcommand{\poly}{polynomial } \newcommand{\polys}{polynomials }
    \newcommand{\pz}{\frac{\p}{\p z}} \newcommand{\pzi}{\frac{\p}{\p z_i}}
    \newcommand{\edge}{\text{\raisebox{2.75pt}{\makebox[20pt][s]{\hrulefill}}}}
    \newcommand{\halfedge}{\text{\raisebox{2.75pt}{\makebox[10pt][s]{\hrulefill}}}}
    \newcommand{\n}{\notag}

\tableofcontents

\renewcommand{\theequation}{\thesection.\arabic{equation}}
\renewcommand{\therema}{\thesection.\arabic{rema}}
\setcounter{equation}{0}
\setcounter{rema}{0}

\vskip 10pt\noindent {\bf Dedication.} This paper is submitted in honor of
Professor Masayoshi Miyanishi, who has profoundly impacted the field
of Affine Algebraic Geometry, giving inspiration to many
mathematicians working in this and related areas, including this
author.

\section{The Jacobian Conjecture}\label{JacConj}
\subsection{The General Assertion} 
The Jacobian Conjecture is:
\begin{conj}[JC]\label{JC} For any integer $n\ge1$ and polynomials
$F_{1},\ldots,F_{n}\in\mathbb C[X_{1},\ldots,X_{n}]$, the polynomial
map $F=(F_{1},\ldots,F_{n}):\mathbb C^{n}\to\mathbb C^{n}$ is an
automorphism if the determinant $|JF|$ of the Jacobian matrix
$JF=(D_{i}F_{j})$ is a nonzero constant.
\end{conj}
Here and throughout this paper we write $D_{i}$ for
${\partial}/{\partial X_{i}}\,$.  We will continue to write $JF$ for
the Jacobian matrix of a polynomial map $F$, and the determinant of
this matrix will be denoted by $|JF|$.

For technical reasons, it will be convenient to henceforth consider
polynomial maps (and later power series maps) with coefficients in an
arbitrary commutative $\mathbb Q$-algebra $K$.  Proving the Jacobian 
Conjecture is equivalent to proving the conjecture as stated in 
\ref{JC} with $\mathbb C$ replaced by $K$ (and $\mathbb C^{n}$ by
$\text{Spec}\,K[X_{1,}\ldots,,X_{n}]$).

\subsection{The Homogeneous Symmetric 
    Reduction}\label{PotSec}

This paper is based on the following result, which puts together two 
well-known reductions:

\begin{theo}[Symmetric Reduction]\label{SymRed} The Jacobian
Conjecture is true if it holds for all polynomial maps $F$ having the
form $F=X-H$ with $H$ homogeneous of degree $d\ge2$ and $JH$ is a symmetric
matix.  In fact, it suffices to prove the case $d=3$.
\end{theo}

\noindent The reduction to the homogeneous cubic case was proved in \cite{BCW};  
the reduction to the symmetric situation is due to de 
Bondt and van den Essen \cite{BE1}

\begin{defi}\label{symhomdefi} A polynomial map $F=X-H$ of the form
prescribed in Theorem \ref{SymRed}, with $d\ge2$ arbitrary, will be
said to be of {\it symmetric homogeneous type.}
\end{defi}

The condition $JH$ is symmetric is equivalent to the existence of a
homogenous polynomial $P\in K[X_{1},\ldots,X_{n}]$ with
$H=\triangledown P$.  $P$ is called the {\it potential function} for
$H$.  Thus the symmetric case occurs precisely when the Jacobian
matrix of $H$ is the Hessian matrix of $P\,$:
    $$JH=\mathrm{Hess}\,P=(D_{i}D_{j}P)$$
If $H$ is homogeneous of degree $d$, $P$ can, of course, be taken to be
homogeneous of degree $d+1$.

\renewcommand{\theequation}{\thesection.\arabic{equation}}
\renewcommand{\therema}{\thesection.\arabic{rema}}
\setcounter{equation}{0}
\setcounter{rema}{0}

\section{Formulas for the Formal Inverse} The formulas for the formal
inverse given in this section provide means for the Jacobian
Conjecture to be addressed as a problem in combinatorics.  See
\cite{W6} for a full discussion of this approach.  These formulas are
valid for systems of power series $F=(F_{1},\ldots,F_{n})$ where, for
$i=1,\ldots,n$, $F_{i}\in K[[X_{1,},\ldots,X_{n}]]$ has the form
$X_{i}+\text{higher degree terms}$.  We call such a map a {\it formal
map of special type.} Such a map has a unique {\it formal inverse,}
that is, a formal map $F^{-1}=G=(G_{1},\ldots,G_{n})$ of special type
having the property that $F\circ G=G\circ F=(X_{1,},\ldots,X_{n})$.

\subsection{The Tree Formula of Bass-Connell-Wright}\label{BCWsection}
Let $\mathbb T_{\text{rt}}$ be the set of isomorphism classes of
finite rooted trees.  For $G=F^{-1}$, the Tree Formula of
Bass-Connell-Wright (best reference for this is \cite{W2}) states:
\begin{theo}[BCW Tree Formula]\label{BCWthm} Let $F=X-H$ be a formal
map, and let $G=(G_{1},\ldots,G_{n})$ be the formal inverse.  Writing
$G=X+N$,
    with $N=(N_{1},\ldots,N_{n})$, we have $$
    N_{i}=\sum_{T\in\mathbb
    T_{\text{rt}}}\frac{1}{|\text{Aut}\,\,T|}\,\,\mathcal P_{T,H,i} $$
    where $$\mathcal
    P_{T,H,i}=\sum_{\substack{\ell:V(T)\to\{1,\ldots,n\} \\
    \ell(\text{rt}_{T})\,=\,i}} \,\,\, \prod_{v\in
    V(T)}D_{\ell(v^{+})}H_{\ell(v)}\,.$$ In this expression $v^{+}$ is
    the set $\{w_{1},\ldots,w_{t}\}$ of children of $v$ and
    $D_{\ell(v^{+})}=D_{\ell(w_{1})}\cdots D_{\ell(w_{t})}\,$.
    \end{theo}
    
    In the case where $F=X-H$ with $H$ homogeneous of degree $d\ge2$,
    the polynomial $\mathcal P_{T,H,i}$ is homogeneous of degree
    $m(d-1)+1$ where $m=|V(T)|$, the number of vertices in the tree
    $T$.  Hence letting $\mathbb T_{{\text{rt}}_{m}}$ be the set of
    trees in $\mathbb T_{\text{rt}}$ having $m$ vertices, and 
    letting  $\mathcal P_{T,H}=(\mathcal P_{T,H,1},\ldots,\mathcal 
    P_{T,H,n})$, we have:
    \begin{theo}[Bass-Connell-Wright Homogeneous Tree
    Formula]\label{BCWthmHomog} Let $F=X-H$ be a polynomial map with 
    $H$ homogeneous of degree $d$.  Then the formal inverse has
    the form $G=X+N$ where
    $$N=N^{(1)}+N^{(2)}+N^{(3)}+\cdots$$
    with $N^{(m)}$ homogeneous of degree $m(d-1)+1$ and given by the formula
    \begin{equation}\label{Nm} N^{(m)}=\sum_{T\in{\mathbb
    T_{\text{rt}}}_{m}}\frac{1}{|\text{Aut}\,\,T|}\,\,\mathcal P_{T,H}\,.
	\end{equation}
    \end{theo} 
    
    \subsection{The Tree Formula for the Symmetric Case} The formula
    of Bass-Connell-Wright takes on a simpler form in the symmetric
    case.  We now let $\mathbb T$ be the set of isomorphism classes of finite
    free trees (i.e., having no designated root).  
    
    \begin{theo}[Symmetric Tree
    Formula]\label{SymTreeThm} Let $F=X-\triangledown P$ be a
    symetric formal map, and let
    $G=(G_{1},\ldots,G_{n})$ be its inverse.  Then $G=X+\triangledown
    Q$ with $$
    Q=\sum_{T\in\mathbb T}\frac{1}{|\text{Aut}\,\,T|}\,\,\mathcal
    Q_{T,P}$$ where \begin{equation}\label{defQ}\mathcal
    Q_{T,P}=\sum_{\ell:E(T)\to\{1,\ldots,n\}} \,\,\, \prod_{v\in
    V(T)}D_{{\rm{adj}}(v)}P\,.\end{equation} Here $\text{adj}(v)$ is
    the set $\{e_{1},\ldots,e_{s}\}$ of edges adjacent to $v$ and
    $D_{\text{adj}(v)}=D_{\ell(e_{1})}\cdots D_{\ell(e_{s})}\,$.
    \end{theo}
    
    \noindent A somewhat similar formula appears without proof
    in \cite{Me}.
    
    \begin{proof}  In the case where $H=\del P$, the expression 
    $H_{l(v)}$ becomes $D_{l(v)}P$, hence 
    \begin{align}
	\mathcal P_{T,H} &=\prod_{v\in V(T)}D_{k(v)}D_{l(v)}P
	\notag\\
	&=\prod_{v\in V(T)}D_{k(v)+e_{l(v)}}P\,.\notag
	\end{align}
    Given $i\in\{1,\ldots,n\}$, $T\in\mathbb T$,
    $\ell:E(T)\to\{1,\ldots,n\}$, and $w\in V(T)$, we create a rooted
    tree $T_{w}$ by declaring $w$ to be the root, and create a
    labeling $\ell_{w}:V(T)\to\{1,\ldots,n\}$ by giving $w$ the label
    $i$ and moving the label of each edge $e\in E(V)$ to the vertex
    $v$ adjacent to $e$ which is farthest from $w$.  Let $k_{w}(v)$ be
    the child type of $v$ in $T_{v}$ resulting from this labeling.
    
    We claim that $Q$ as defined in the theorem is the 
    potential function for $$N=\sum_{T\in\mathbb
     T_{\text{rt}}}\frac{1}{|\text{Aut}\,\,T|}\,\,\mathcal P_{T,H}\,,$$  
    that is to say, $D_{i}Q=N_{i}$ for $i=1,\ldots,n$.  To see this, 
    note that: 
    \allowdisplaybreaks{
    \begin{align}
    D_{i}Q&= D_{i}\sum_{T\in\mathbb
    T}\frac{1}{|\text{Aut}\,\,T|}\,\,\mathcal Q_{T,P}\notag\\
    &=D_{i}\left(\sum_{T\in\mathbb
    T}\frac{1}{|\text{Aut}\,\,T|}\,\,\sum_{\ell:E(T)\to\{1,\ldots,n\}}
    \,\,\, \prod_{v\in V(T)}D_{\text{adj}(v)}P\right)\notag\\
    &=\sum_{T\in\mathbb
    T}\frac{1}{|\text{Aut}\,\,T|}\,\,\sum_{\ell:E(T)\to\{1,\ldots,n\}}
    \,\,\, \sum_{w\in V(T)}\,\prod_{v\in 
    V(T)}D_{\text{adj}(v)+\delta_{v,w}e_{i}}P\notag\\
    &=\sum_{T\in\mathbb
	T}\frac{1}{|\text{Aut}\,\,T|}\,\,\sum_{\ell:E(T)\to\{1,\ldots,n\}}
	\,\,\, \sum_{w\in V(T)}\,\prod_{v\in 
	V(T)}D_{k_{w}(v)}D_{\ell_{w}(v)}P\notag\\
	&=\sum_{T\in\mathbb
		T}\frac{1}{|\text{Aut}\,\,T|}\,\,\sum_{\ell:E(T)\to\{1,\ldots,n\}}
		\,\,\, \sum_{w\in V(T)}\,\prod_{v\in 
		V(T)}D_{k_{w}(v)}(\del P)_{\ell_{w}(v)}\notag\\
			&=\sum_{T\in\mathbb
			T}\frac{1}{|\text{Aut}\,\,T|}
			\,\,\sum_{w\in 
			V(T)}\,\,\,\sum_{\substack{h:V(T)\to\{1,\ldots,n\}\\h(w)=i}}\,\prod_{v\in 
		V(T)}D_{k_{w}(v)}(\del P)_{h(v)}\notag\\
		&=\sum_{T\in\mathbb
		T}\frac{1}{|\text{Aut}\,\,T|}
		\,\,\sum_{w\in
		V(T)}\,\,\,\sum_{\substack{h:V(T)\to\{1,\ldots,n\}\\h(w)=i}}\,\mathcal 
		P_{T_{w},\del P,i}\notag\\
		&=\sum_{S\in \mathbb T_{\text{rt}}}\,\,\sum_{T\in\mathbb 
		T}\,\sum_{\substack{w\in V(T)\\T_{w}\cong 
		S}}\,\frac{1}{|\text{Aut}\,\,T|}\,\,\mathcal 
		P_{S,\del P,i}\notag\\
		\intertext{Denoting by $\bar S$, for $S\in \mathbb T_{\text{rt}}$, 
		the unrooted tree determined by $S$, ignoring the root, we have}
		&=\sum_{S\in \mathbb T_{\text{rt}}}\,\,\sum_{\substack{w\in V(\bar 
		S)\\\bar S_{w}\cong_{\mathbb 
		T_{\text{rt}}}S}}\,\frac{1}{|\text{Aut}\,\,\bar S|}\,\,\mathcal 
		P_{S,\del P,i}\notag\\
		&=\sum_{S\in \mathbb T_{\text{rt}}}\,\frac{\left|\{w\in V(\bar 
		S)\,|\, \bar S_{w}\cong_{\mathbb 
		T_{\text{rt}}}S\}\right|}{|\text{Aut}\,\,\bar S|}\,\,\mathcal 
		P_{S,\del P,i}\notag\\
		\intertext{Aut$\,\bar S$ acts on $V(\bar S)=V(S)$, the
		orbit of the root $r$ of $S$ being the set $\{w\in
		V(\bar S)\,|\, \bar S_{w}\cong_{\mathbb
		T_{\text{rt}}}S\}$.  The stabilizer of $r$ in
		Aut$\,\bar S$ is Aut$_{\mathbb T_{\text{rt}}}\,S$, so
		$$\left|\{w\in V(\bar S)\,\right|\, \bar
		S_{w}\cong_{\mathbb
		T_{\text{rt}}}S\}|=\frac{|\text{Aut}\,\bar S|}{|\text{Aut}_{\mathbb 
		T_{\text{rt}}}\,S|}$$ and we get }
		&=\sum_{S\in \mathbb T_{\text{rt}}}\,\,\frac{1}{|\text{Aut}_{\mathbb 
		T_{\text{rt}}}\,S|}\,\,\mathcal 
		P_{S,\del P,i}\notag\\&=N_{i}\,,\notag
	\end{align}
	}
	which, since $\del P=H$, completes the proof.
    \end{proof}
    
    For the symmetric homogeneous case Theorem \ref{SymTreeThm} gives
    the following.  Here we let $\mathbb T_{m}$ be the set of
    isomorphism classes of free (i.e.,
    non-rooted) trees having $m$ vertices.
    \begin{theo}[Symmetric Homogeneous Tree Formula]\label{SymHomThm}
    Suppose $F$ has the form $F=X-\del P$ with $P$ homogeneous of
    degree $d+1$.  Let $G$ be the formal inverse of $F$.  Then
    $G=X+\triangledown Q$ with $$Q=Q^{(1)}+Q^{(2)}+Q^{(3)}+\cdots$$ 
    and $$Q^{(m)}=\sum_{T\in\mathbb
    T_{m}}\frac{1}{|\text{Aut}\,\,T|}\,\,\mathcal Q_{T}\,.$$
    $Q^{(m)}$ is homogeneous of degree $m(d-1)+2$. \end{theo}
    
    \noindent It is clear that in this situation
    $N^{(m)}=\triangledown Q^{(m)}$, where $N^{(m)}$ is as in Theorem 
    \ref{BCWthm}.

    \subsection{Zhao's Formulas and the Gap Theorem} The formula below
    of Zhao, proved in \cite{Z3}, has an important consequence for
    this discussion, namely the Gap Theorem (Theorem \ref{gapthm}).
    
	\begin{theo}[Zhao's Formula for the Symmetric Case]\label{Zhaothm2} As
	in Theorem \ref{SymHomThm}, let $Q^{(m)}$, $m\ge1$, be the homogeneous
	summands of the potential function for $N=G-X$, where $G$ is formal
	inverse of a degree $d$ polynomial map $F=X-\triangledown P$
	of symmetric homogeneous type.  Then $Q^{(1)}=P$ and, for $m\ge2$,
	\begin{equation}\label{zform2}Q^{(m)}=\frac{1}{2(m-1)}
	\sum_{\substack{k+\ell=m\\k,\ell\ge1}}\left(\triangledown Q^{(k)}\cdot
	\triangledown Q^{(\ell)}\right)\,.\end{equation} (Here $(\triangledown
	Q^{(k)}\cdot \triangledown Q^{(\ell)})$ denotes the usual dot product
	of vectors.)  \end{theo} \noindent Again it should be noted that this
	theorem holds in the nonhomogeneous case as well, giving
	nonhomogeneous, formally converging summands for the potential function for $N$.
	
	The following theorem gives explicit finitude to showing that the polynomial 
	inverse of a polynomial map of symmetric homogeneous type is a 
	polynomial.
	
	\begin{theo} [Gap Theorem for the Symmetric Case] \label{gapthm} Given the
	situation of Theorem \ref{SymHomThm}, then $F$ is invertible,
	i.e., $G$ is a polynomial map, if we have
	$$Q^{(M+1)}=Q^{(M+2)}=\cdots=Q^{(2M)}=0$$
	for some positive integer $M$.
	    \end{theo}
	    \begin{proof}  This is immediate from formula \ref{zform2} in 
	    Theorem \ref{Zhaothm2}.
		\end{proof}

\renewcommand{\theequation}{\thesection.\arabic{equation}}
\renewcommand{\therema}{\thesection.\arabic{rema}}
\setcounter{equation}{0}
\setcounter{rema}{0}

\section{Consequences}\label{IMQ}  

\subsection{Trees with naked chains} In the case where $F=X-H$ with
$H$ homogeneous of degree $d\ge2$, then the invertibility of $JF$ is
equivalent to $JH$ being nilpotent, in which case we must have
$(JH)^{n}=0$ (where $n$ is the number of variables).  This motivates
the following theorem:

\begin{theo}[Chain Vanishing Theorem] \label{ChainVanThm}  Suppose $P\in 
K[[X_{1},\ldots,X_{n}]]$ with $(\text{Hess}\,P)^{r}=0$ for some 
$r\ge1$, 
and suppose $T$ is a tree which contains a 
``naked $r$-chain,'' that is, a geodesic 
$$(\overset{e_{0}}\edge)\underset{v_{1}}\bullet\overset{e_{1}}\edge\underset{v_{2}}\bullet\overset{e_{2}}\edge
\underset{v_{3}}\bullet\halfedge\,\cdots\,\halfedge
\underset{v_{r-1}}\bullet\overset{e_{r-1}}\edge\underset{v_{r}}\bullet(\overset{e_{r}}\edge)$$ 
meaning the vertices $v_{2}.\ldots,v_{r-1}$ have degree 2 and the two 
vertices $v_{1},v_{r}$ have degree 1 or 2. Assume either (a) $P$ is a 
homogeneous polynomial of degree $\ge2$, or (b) both 
$v_{1}$ and $v_{r}$ have degree 2.  Then $\mathcal Q_{T,P}=0$.
\end{theo}
\proof First we assume (b) holds, i.e., $e_{0}$ and $e_{r}$ are 
actually there.  Write $E(T)$ as the disjoint union 
$\{e_{1},\ldots,e_{r-1}\}\cup E'$ and $V(T)$ as the disjoint union $\{v_{1},\ldots,v_{r}\}\cup 
V'$.  By definition (see Theorem \ref{SymTreeThm}) we have
{\allowdisplaybreaks
\begin{align}\mathcal Q_{T,P}&=\sum_{\ell:E(T)\to\{1,\ldots,n\}}
    \,\,\, \prod_{v\in V(T)}D_{\text{adj}(v)}P\n\\
    &=\sum_{\ell':E'\to\{1,\ldots,n\}}
    \,\,\,\sum_{\ell:\{e_{1},\ldots,v_{r-1}\}\to\{1,\ldots,n\}} \,\,\,
    \prod_{v\in V'}D_{\text{adj}(v)}P\,\,\, \prod_{v\in
    \{v_{1},\ldots,v_{r}\}}D_{\text{adj}(v)}P\n\\ 
    &=\sum_{\ell':E'\to\{1,\ldots,n\}} \,\,\,\prod_{v\in
    V'}D_{\text{adj}(v)}P\,\,\,\sum_{\ell:\{e_{1},\ldots,v_{r-1}\}\to\{1,\ldots,n\}}
    \,\,\,\prod_{v\in
    \{v_{1},\ldots,v_{r}\}}D_{\text{adj}(v)}P\n\\
    &=\sum_{\ell':E'\to\{1,\ldots,n\}} \,\,\,\prod_{v\in
	V'}D_{\text{adj}(v)}P\n\\
	&\hskip2cm\sum_{i_{1},\ldots,i_{r-1}}
	\left(D_{\ell'(e_{0})\,i_{1}}P\right)\left(D_{i_{1}i_{2}}P\right)
	\cdots\left(D_{i_{r-2}i_{r-1}}P\right)
	\left(D_{i_{r-1}\ell'(e_{r})}P\right)\label{prod}
\end{align}
}
Since the $(ij)^{\text{th}}$ entry in the matrix $\text{Hess}\,P$ is 
$D_{ij}P$, the final summation above gives the 
$(\ell'(e_{0})\,\ell'(e_{r}))^{\text{th}}$ entry in 
$(\text{Hess}\,P)^{r}$, which is zero by hypothesis.  Therefore $\mathcal 
Q_{T,P}=0$.

Now assume (a) holds.  We proceed as before and all the equalities 
above are valid except the last one, which assumes the existence of 
$e_{0}$ and $e_{r}$.  If, say, $e_{r}$ is present but $e_{0}$ is not, 
then the final summation \ref{prod} reads:
\begin{equation}\sum_{i_{1},\ldots,i_{r-1}}
	\left(D_{i_{1}}P\right)\left(D_{i_{1}i_{2}}P\right)
	\cdots\left(D_{i_{r-2}i_{r-1}}P\right)
	\left(D_{i_{r-1}\ell'(e_{r})}P\right)\,.\n\end{equation}
	Since $D_{i_{1}}P$ is homogeneous of degree $d-1$, Eulers formula 
	says $D_{i_{1}}P=\frac{1}{d-1}\sum_{i_{0}=1}^{n}D_{i_{0}i_{1}}P$.  
	Thus the above sum is
	$$\frac{1}{d-1}\sum_{i_{0},i_{1},\ldots,i_{r-1}}
		\left(D_{i_{0}i_{1}}P\right)\left(D_{i_{1}i_{2}}P\right)
		\cdots\left(D_{i_{r-2}i_{r-1}}P\right)
		\left(D_{i_{r-1}\ell'(e_{r})}P\right)\,,$$
		which vanishes, since $(\text{Hess}\,P)^{r}=0$.  Finally, of both 
		$e_{0}$ and $e_{r}$ are absent, then \ref{prod} becomes
		\begin{equation}\sum_{i_{1},\ldots,i_{r-1}}
			\left(D_{i_{1}}P\right)\left(D_{i_{1}i_{2}}P\right)
			\cdots\left(D_{i_{r-2}i_{r-1}}P\right)
			\left(D_{i_{r-1}}P\right)\,,\n\end{equation}
			and the proof is completed by applying Euler's formula to both end 
			factors $D_{i_{1}}P$ and $D_{i_{r}}P$.
\endproof
	
\subsection{The Symmetric $JH^{3}=0$ Case} 
	
The following new result for the
symmetric situation, announced in \cite{W6}, will
use the Symmetric Homogeneous Tree Formula and the Chain Vanishing Theorem 
(Theorems \ref{SymHomThm} and \ref{ChainVanThm}).

\begin{theo}[Symmetric Cube Zero Case] \label{JH3=0Sym} If $F=X-H$ is a
polynomial map with symmetric Jacobian matrix of homogeneous type with
$(JH)^{3}=0$, then $F$ is invertible with $$F^{-1}=X+N^{(1)}
+N^{(2)}\,.$$ In particular, the degree of $F^{-1}$ is $\le2d-1$, 
where $d=\text{deg}\,H$ 
(independent of $n$).
\end{theo}

\begin{rema} \label{rem1} What is remarkable about the above statement is that it is 
independent of $n$, the number of variables.  Moreover the form of 
$F^{-1}$ is independent of the degree $d$ of $H$.  (The known bound 
for the degree of the inverse of an invertible polynomial map of 
degree $d$ is $d^{n-1}$ (Gabber's Theorem).  See \cite{BCW}.)\end{rema}

\begin{proof} By Gap Theorem (Theorem \ref{gapthm}) it sufffices to
show that $Q^{(3)}=Q^{(4)}=0$, where $Q^{(m)}$ is as defined in the
Symmetric Homogeneous Tree Formula (Theorem \ref{SymHomThm}).  But this is
immediate from the following proposition.
    \end{proof}
    
    \begin{propo} \label{JH3=0TreeProp}  If $P\in K[X_{1},\ldots,X_{n}]$ is homogeneous of 
    degree $\ge2$ with $(\text{Hess}\,P)^{3}=0$ and if $T$ is a tree 
    with 3 or 4 vertices, then $\mathcal Q_{T,P}=0$.
	\end{propo}
	
	\begin{proof}  The only tree with three vertices is the 
	3-chain $T=\bullet\edge\bullet\edge\bullet$, and in this case 
	$\mathcal Q_{T,P}=0$ by the Chain Vanishing theorem 
	(\ref{ChainVanThm}).  There are two trees with four vertices, namely
	\begin{equation}T_{1}=\bullet\edge\bullet\edge\bullet\edge\bullet\qquad\text{and}\qquad
	    T_{2}=\bullet\edge\overset{\overset{\displaystyle{\bullet}}
	    {\displaystyle{|}}}\bullet\edge\bullet\n\end{equation}
	    We have $\mathcal Q_{T_{1},P}=0$ by the Chain Vanishing Theorem.  
	    To get the vanishing of $\mathcal Q_{T_{2},P}$ we apply the 
	    operator $\sum_{i=1}^{n}(D_{i}P)D_{i}$ to $\mathcal Q_{T,P}(=0)$, where 
	    $T$, as above, is the 3-chain.  We get:
	    \begin{align}
		0&=\sum_{i=1}^{n}(D_{i}P)(D_{i}\mathcal Q_{T,P})\n\\
		&=\sum_{i=1}^{n}(D_{i}P)\left(D_{i}\sum_{j,k}(D_{j}P)(D_{jk}P)(D_{k}P)\right)\n\\
		\intertext{which becomes, using the product rule:}
		&=\sum_{i,j,k}(D_{i}P)(D_{ij}P)(D_{jk}P)(D_{k}P)\label{1st}\\
		&\qquad\qquad+\sum_{i,j,k}(D_{i}P)(D_{j}P)(D_{ijk}P)(D_{k}P)\label{2nd}\\
		&\qquad\qquad\qquad\qquad+\sum_{i,j,k}(D_{j}P)(D_{jk}P)(D_{ki}P)(D_{i}P)\,.\label{3rd}
		\end{align}
		Note that \ref{1st} and \ref{3rd} are each equal to $\mathcal 
		Q_{T_{1},P}$ and the \ref{2nd} is $\mathcal Q_{T_{2},P}$.  Thus we 
		have $2\mathcal Q_{T_{1},P}+\mathcal Q_{T_{2},P}=0$.  Since $\mathcal 
		Q_{T_{1},P}=0$, we must have $\mathcal Q_{T_{2},P}=0$ as well.
	\end{proof}

\subsection{The Grossman-Larson Algebra}

The proof of Proposition \ref{JH3=0TreeProp} entails operations that 
hearken to a ring defined by Grossman and Larson in 
\cite{GL}, which we will now define as a $\mathbb Q$-algebra.  

Let $\mathcal H_{GL}$, or simply $\mathcal H$, be the vector space
over $\mathbb Q$ spanned by $\mathbb T_{\text{rt}}$, the set of all
rooted trees.  To explain multiplication in $\mathcal H$ it will
be necessary to introduce some concepts and notations.

First, let $S$ be a rooted tree, $T$ a (possibly non-rooted) tree, and
let $v\in V(T)$.  We denote by $S\multimap_{v}T$ the tree which joins
$T$ to $S$ by introducing a new edge $e$ which connects
$\text{rt}_{S}$ to $v$.  If $T$ is a rooted tree, then
$S\multimap_{v}T$ is rooted by $\text{rt}_{T}$.  Similarly if
$S_{1},\ldots,S_{r}$ are rooted trees and $v_{1},\ldots,v_{r}\in
V(T)$, we can form the tree
$$(S_{1},\ldots,S_{r})\multimap_{(v_{1},\ldots,v_{r})}T\,,$$ which
attaches $S_{i}$ to $T$ at $v_{i}$, for $i=1,\ldots,r$.  Again, if $T$
is rooted, we take $\text{rt}_{T}$ to be the root of the newly formed
tree.  

Secondly, if $S$ is a rooted tree, let $\text{DelRoot}(S)$ denote the
{\it forest} (meaning a set with multiplicity of rooted trees) of branches
of $\text{rt}_{S}$.  This means we delete the root of $S$ and its
adjacent edges; the children of $\text{rt}_{S}$ become the roots of 
the trees in $\text{DelRoot}(S)$.

Now we define the multiplcation in $\mathcal H$.  For rooted 
trees $S$ and $T$, we write $\text{DelRoot}(S)=\{S_{1},\ldots,S_{r}\}$ 
(incorporating multiplicity) and we define the product $S\cdot T$ by
\begin{equation}\label{mult}
    S\cdot T=\sum_{(v_{1},\ldots,v_{r})\in V(T)^{r}}
    \left[(S_{1},\ldots,S_{r})\multimap_{(v_{1},\ldots,v_{r})}T\right]
    \end{equation}
    This multiplication is extended to $\mathcal H$ by distributivity.  One
    quickly checks that the singleton serves as a left and right
    multiplicative identity element.  In \cite{GL} it is shown that
    the multiplication is associative (a fact which is not hard to
    verify), and that $\mathcal H$ has the additional structure of a
    Hopf algebra, a property which will not be used here.
    
    An important thing to note is that $\mathcal H$ is a graded ring 
    by the grading
    $\mathcal H=\oplus_{i=0}^{\infty}\,\mathcal H_{i}$, where
    $\mathcal H_{i}$ spanned by trees having $i$ non-root vertices, 
    i.e., by $\mathbb T_{\text{rt}_{i+1}}$.
    
    Now we let $\mathcal M$ be the $\mathbb Q$-vector space spanned by
    the set $\mathbb T$ of all non-rooted trees.  We observed that
    $S\multimap_{v}T$ forms a non-rooted tree when $S$ is rooted and
    $T$ is non-rooted; one can use \ref{mult} to endow $\mathcal M$
    with the structure of an $\mathcal H$-module, which we will call 
    the {\it tree module.}.  In fact, $\mathcal
    M$ is a graded $\mathcal H$-module $\mathcal
    M=\oplus_{i=1}^{\infty}\,\mathcal M_{i}$ taking $\mathcal M_{i}$ 
    to be the vector space spanned by $\mathbb T_{i}$.
    
    \begin{defi} [Free tree quotient modules] \label{submod} For a positive integer $r$, let
    $\mathcal C(r)$ denote the sub-$\mathcal H$-module of $\mathcal M$
    generated by all trees containing a naked $r$-chain (see Theorem
    \ref{ChainVanThm} for the definition).  Let $\mathcal V(r)$ denote
    the sub-vector space (over $\mathbb Q$) generated by all trees
    which have at least one vertex of degree $\ge r+1$.  It is easily
    seen that $\mathcal V(r)$ is also a sub-$\mathcal H$-module.  For
    positive integers $r,e$, let $\mathcal N(r,e)=\mathcal
    C(r)+\mathcal V(e)$.  These are graded submodules of $\mathcal M$.
    Finally, let $\overline{\mathcal M}(r,e)=\mathcal M/\mathcal
    N(r,e)$ and let $\overline{\mathcal M}(r,\infty)=\mathcal
    M/\mathcal C(r)$.  The $\mathcal H$-modules $\overline{\mathcal 
    M}(r,e)$ ($e$ possibly being $\infty$) will be called the {\it 
    tree quotient modules.} \end{defi}
    
    Given $\gamma\in\mathcal M$ we will often denote by
       $\overline\gamma$ its image in $\overline{\mathcal M}(r,e)$, where
       $r$ and $e$ are understood in the context of the discussion.
    
    \subsection{Relationship to the Ring of Differential Operators}
    
    We write $\mathfrak D[X]=\mathfrak D[X_{1},\ldots,X_{n}]$ for the
    ring of differential operators on $K[X]=K[X_{1},\ldots,X_{n}]$.  A
    polynomial $P\in K[X]$ gives rise to a ring homomorphism
    \begin{equation} \varphi_{P}:\mathcal H\to \mathfrak
    D[X]\end{equation} which we will be defined as follows: For a
    rooted tree $S$, we let $e_{1},\ldots, e_{r}$ be the edges 
    adjacent to $\text{rt}_{S}$ and define the differential operator $\mathfrak 
    d_{S,P}\in 
    \mathfrak D[X]$
    by $$\mathfrak d_{S,P}=\sum_{\ell:E(S)\to\{1,\ldots,n\}} \,\,\,
    \left(\prod_{v\in
    V(S)-\{\text{rt}_{S}\}}D_{\text{adj}(v)}P\right)D_{\ell(e_{1})\ell
    (e_{2})\cdots\ell(e_{r})}\,.$$ Note the similarity with the 
    definition of the polynomial $\mathcal
    Q_{T,P}$ (Theorem \ref{SymTreeThm}) for a free tree $T$; the difference is that here we
    omit $\text{rt}_{S}$ from the product and leave ``open'' the derivatives 
    corresponding to edges adjacent to $\text{rt}_{S}$.
    
    Taking $\varphi_{P}(S)=\mathfrak d_{S,P}$ defines $\varphi_{P}$ on
    $\mathcal H$ as a $\mathbb Q$-linear map; in fact, it is straightforward 
    to show that $\varphi_{P}$ is a ring homomorphism.
    
    Now we define a map \begin{equation}\rho_{P}:\mathcal M\to 
    K[X]\end{equation} by sending an unrooted tree $T$ to $\mathcal 
    Q_{T,P}$.  Again, it is straightforward to verify that this map is 
    compatible with the structures of $\mathcal M$ as an $\mathcal 
    H$-module and $K[X]$ as a $\mathfrak D[X]$-module, that is, the 
    diagram
    \begin{equation}
	\begin{matrix}\mathcal H\times\mathcal M & \to & \mathcal M\\
	    \downarrow & {} &\downarrow\\
	    \mathfrak D[X]\times K[X] & \to & K[X]
	    \end{matrix}
	    \end{equation}
commutes, where the horizontal arrows are induced by the module 
structures and the vertical arrows are $\varphi_{P}\times\rho_{P}$ 
and $\rho_{P}$.  Now we observe:

\begin{propo}\label{rhobar} For $P\in K[X_{1},\ldots,X_{n}]$ and positive integers 
$r,e$ we have:
\begin{enumerate}
    
    \item If $P$ is homogeneous with $\text{Hess}\,(P)^{r}=0$, then 
    $\rho_{P}(\mathcal C(r))=0$.
    
    \item If $\text{deg }P\le e$, then $\rho_{P}(\mathcal V(e))=0$.
\end{enumerate}
Thus if $P$ is homogeneous of degree $\le e$ with
    $\text{Hess}\,(P)^{r}=0$ then $\rho_{P}$ induces a homomorphism 
    $\overline{\rho}_{P}(r,e):\overline{\mathcal M}(r,e)\to K[X]$ such 
    that
    \begin{equation}\label{defrhobar}
	    \begin{matrix}\mathcal H\times\overline{\mathcal M}(r,e) &
	    \to & \overline{\mathcal M}(r,e)\\
		\downarrow & {} &\downarrow\\
		\mathfrak D[X]\times K[X] & \to & K[X]
		\end{matrix}
		\end{equation}
commutes, where the horizontal arrows are induced by the module 
structures and the vertical arrows are $\varphi_{P}\times\overline{\rho}_{P}(r,e)$
and $\overline{\rho}_{P}(r,e)$. The last statement also holds for 
$e=\infty$.
    \end{propo}
    
\begin{proof} Statement 1 follows from the Chain Vanishing Theorem 
(\ref{ChainVanThm}).  Statement 2 follows from the definition of 
$\mathcal Q_{T,P}$ (\ref{defQ}) with the observation that 
$D_{\text{adj}(v)}P=0$ if $v\in V(T)$ has degree $\ge e+1$.
    \end{proof}
    
\begin{defi}\label{defnu} For $m\ge1$ let $\nu_{m}\in\mathcal M$ be defined by 
$$\nu_{m}=\sum_{T\in\mathbb 
T_{m}}\frac{1}{\left|\text{Aut}\,T\right|}\,T\,.$$  
    \end{defi}
    
    Note that $\nu_{m}$ is homogeneous of degree $m$, i.e.,
    $\nu_{m}\in\mathcal M_{m}$, and that, for $P\in K[X]$ homogeneous,
    $\rho_{P}(\nu_{m})=Q^{(m)}$, where $Q^{(m)}$ is as defined in
    Theorem \ref{SymTreeThm}.  It follows that
    \begin{equation}\label{rhonu}\overline{\rho}_{P}(\overline{\nu}_{m})=Q^{(m)}\,\end{equation}
    where $\overline{\nu}_{m}$ is the image of $\nu_{m}$ in
    $\overline{\mathcal M}(r,e)$, whenever $\overline{\rho}_{P}$ makes
    sence by virtue of Proposition \ref{rhobar}.

\subsection{The Symmetric $JH^{3}=0$ Case Revisited}\label{cube} 

We will now observe that the proof of Theorem \ref{JH3=0Sym} boils
down to a statement about the $\mathcal H$-module $\overline{\mathcal
M}(3,\infty)$.  The theorem followed from the fact that
$Q^{(3)}=Q^{(4)}=0$ when $P$ is homogeneous and
$\text{Hess}(P)^{3}=0$.  Since $\overline{\rho}_{P}$ is defined for
$r=3$, $e=\infty$, in this situation (by Proposition \ref{rhobar}),
this would follow from $\overline{\nu}_{3} = \overline{\nu}_{4} = 0$
in $\overline{\mathcal M}(3,\infty)$, by \ref{rhonu}.  But in fact we
have, more strongly:

\begin{propo}\label{3infty}  In the graded module $\overline{\mathcal 
M}(3,\infty)=\oplus_{i=1}^{\infty}\overline{\mathcal 
M}(3,\infty)_{i}$ the homogeneous summands $\overline{\mathcal 
M}(3,\infty)_{3}$ and $\overline{\mathcal 
M}(3,\infty)_{4}$ are both zero.
    \end{propo}
    
    \begin{proof} Let $T$, $T_{1}$, and $T_{2}$ be as defined in the
    proof of Proposition \ref{JH3=0TreeProp}.  We will write $\bar S$
    for the image in $\overline{\mathcal M}(3,\infty)$ of a tree $S$.
    Since $\mathcal M_{3}=\mathbb Q\cdot T$ and $T$ is the chain of
    length 3, $\bar T=0$ and hence $\overline{M}(3,\infty)_{3}=0$.  We
    have $\mathcal M_{4}=\mathbb Q\cdot T_{1}\,\oplus\,\mathbb Q\cdot
    T_{2}$ and since $\bar T_{1}=0$ (since $T_{1}$ is the chain of
    length 4), $\overline{\mathcal M}(3,\infty)$ is generated over
    $\mathbb Q$ by $\bar T_{2}$.  Now note that, letting $S$ be the
    rooted chain of length 2, i.e.,
    $$S=\bullet\edge\overset{\text{rt}}\bullet$$
    then the $\mathcal H$ action on $\mathcal M$ gives $S\cdot
    T=2T_{1}+T_{2}$, from which it follows that $T_{2}\in\mathcal
    C(3)$.  Therefore $\bar T_{2}=0$ and so
    $\overline{M}(3,\infty)_{4}=0$, completing the proof.
	
	\end{proof}

\subsection{The Quadratic Symmetric $JH^{4}=0$ Case}\label{quad}

Computations in the $\mathcal H$-modules $\overline{\mathcal M}(r,e)$
allow us to obtain certain specific results when $JH$ is nilpotent of
higher order, for certain specific degrees.  For example:

\begin{theo} \label{d=2,JH4=0Sy} Let $F=X-H$ be a polynomial map having symmetric
Jacobian matrix, with $H$ quadratic homogeneous and $(JH)^{4}=0$.
Then $F$ is invertible with
$$F^{-1}=X+N^{(1)}+N^{(2)}+N^{(3)}+N^{(4)}\,.$$
In particular, the degree of $F^{-1}$ is $\le5$.
\end{theo}

\begin{rema} Of course the Jacobian Conjecture is known to be true for
quadratic maps.  This was proved by S. Wang; a simple proof due to S. 
Oda can be found in in \cite{BCW}.  However, Theorem
\ref{d=2,JH4=0Sy} yields more strongly the uniform degree bound of 5 for
$F^{-1}$, when $F$ is as in the theorem, independent of the
number of variables.  Again recall that the general known degree 
bound here is $2^{n-1}$ (see remark \ref{rem1}).\end{rema}

    \begin{proof} The proof will entail an explicit computation in the
    $\mathcal H$-module $\overline{\mathcal M}(4,3)$.  It follows from
    the Gap Theorem (Theorem \ref{gapthm}) that it suffices to prove
    $Q^{(m)}=0$ for $5\le m\le8$, where $Q^{(m)}$ is defined as in
    Theorem \ref{SymHomThm}.
    
    We have $H=\del P$ where $P\in K[X]$ is homogeneous cubic.
    According to Proposition \ref{rhobar} the map
    $\overline{\rho}_{P}(4,3):\overline{\mathcal M}(4,3)\to K[X]$ is 
    defined, with 
    $\overline{\rho}_{P}(4,3)(\overline{\nu}_{m})=Q^{(m)}$, so it 
    suffices to show $\overline{\nu}_{m})=0$ for $5\le m\le8$.  But 
    since $\overline{\nu}_{m}\in\overline{\mathcal M}_{m}(4,3)$
    this follows from the proposition below.
\end{proof}
      
\begin{propo}\label{43}  In the graded module $\overline{\mathcal 
M}(4,3)=\oplus_{i=1}^{\infty}\overline{\mathcal 
M}(4,3)_{i}$ we have $\overline{\mathcal 
M}(4,3)_{m}=0$ for $5\le m\le8$.
    \end{propo}
    
    \begin{proof} Let $S,S',S'',S''',S''''$ and $S'''''$ be the rooted trees
    appearing in Figure \ref{Fi:treesmisc}, the bottom vertex being
    the root.
    \begin{figure}[h]
		    $$\overset{\text{\includegraphics{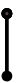}}}{S}\quad
		    \overset{\text{\includegraphics{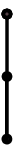}}}{S'}\quad
		    \overset{\text{\includegraphics{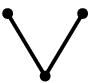}}}{S''}\quad
		    \overset{\text{\includegraphics{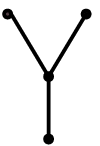}}}{S'''}\quad
		    \overset{\text{\includegraphics{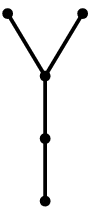}}}{S''''}\quad
		    \overset{\text{\includegraphics{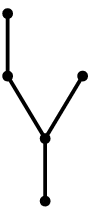}}}{S'''''}
		    $$
		    \caption{}
		    \label{Fi:treesmisc}
		    \end{figure} These will be viewed as elements of the Grossman-Larson algebra $\mathcal 
H$.		    
		    
    Figures \ref{Fi:trees5} through \ref{Fi:trees8} give complete
    lists\footnote{In these lists the trees are ordered by their
    maximal rooted planar representative.  Rooted planar trees are
    ordered by considering first the number of vertices, then, if those
    are the same the lexicographical ordering of the root branches,
    considered left to right.  Inductively, this gives a total
    ordering of rooted planar trees.  Letting $T_{n}$ be the number of rooted trees 
    having $n$ vertices, let $T(x)=\sum_{p=1}^{\infty}T_{p}x^{p}$ be 
    the generating function for rooted trees.  Then $T_{n}$ can be calculated using the 
    following formula, due to G. P\'olya: $$T(x)=
    x\,\text{exp}
	\left\{\sum_{k=1}^{\infty}\frac{T(x^{k})}{k}\right\}$$ Then the 
    number $t_{n}$ of free trees with $n$ vertices is determined 
    by the formula of R. Otter:
    $$t(x)=T(x)-\frac{1}{2}\left\{[T(x)]^{2}-T(x^{2})\right\}$$
    where $t(x)=\sum_{p=1}^{\infty}t_{p}x^{p}$.  The first few values 
    of $t_{n}$ have been 
    found to be: 
    $$t_{1}=1,\,t_{2}=1,\,t_{3}=1,\,t_{4}=2,\,t_{5}=3,\,t_{6}=6,\,t_{7}=11,\,
    t_{8}=23,\,
    t_{9}=47,\,t_{10}=106,\,t_{11}=235$$  This confirms that the lists in 
    Figures \ref{Fi:trees5} through \ref{Fi:trees8} are complete.  See 
    \cite{HP} as a reference for the facts in this footnote.} 
    of free (i.e., unrooted) trees with $m$ vertices,
	\begin{figure}[h]
	    $$\overset{\text{\includegraphics{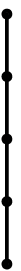}}}{A_{1}}\qquad
	    \overset{\text{\includegraphics{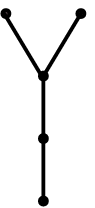}}}{A_{2}}\quad
	    \overset{\text{\includegraphics{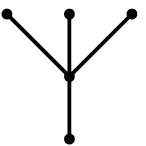}}}{A_{3}}$$
	    \caption{Trees with 5 vertices}
	    \label{Fi:trees5}
	    \end{figure}
	    for $5\le m\le8$.
    Viewing these free trees as elements of $\mathcal M$, our goal is to
    show that each lies in the $\mathcal H$-submodule $\mathcal
    N(4,3)=\mathcal C(4)+\mathcal V(3)$ (see Definition \ref{submod}).  
    
    We first consider the three trees with 5 vertices, identified in
    Figure \ref{Fi:trees5}.  Obviously $A_{1}\in\mathcal C(4)$ and
    $A_{3}\in\mathcal V(3)$.  Furthermore, letting $A$ be the chain
    with four vertices (hence $A\in\mathcal C(4)$), we have $S\cdot
    A=2A_{1}+2A_{2}$, which shows $A_{2}\in\mathcal C(4)$.  Therefore
    $\overline{\mathcal M}(4,3)_{5}=0$.

Figure \ref{Fi:trees6} lists and labels the six trees with 6 vertices.
In $\mathcal M_{6}$,
\begin{figure}[!h]
	    $$\overset{\text{\includegraphics{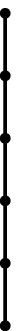}}}{B_{1}}\quad
	    \overset{\text{\includegraphics{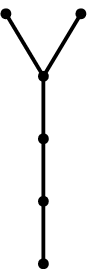}}}{B_{2}}\quad
 	    \overset{\text{\includegraphics{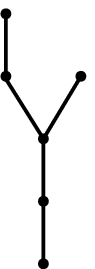}}}{B_{3}}\quad
	    \overset{\text{\includegraphics{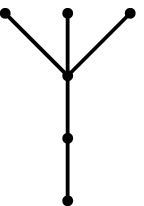}}}{B_{4}}\quad
	    \overset{\text{\includegraphics{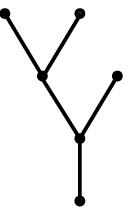}}}{B_{5}}\quad
	    \overset{\text{\includegraphics{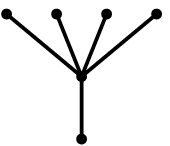}}}{B_{6}}
	    $$
	    \caption{Trees with 6 vertices}
	    \label{Fi:trees6}
	    \end{figure}
	    note that $B_{1}\in\mathcal C(4)$ and that
$B_{4},B_{6}\in\mathcal V(3)$.  Furthermore we have $S'\cdot 
A=2B_{1}+2B_{3}$ which gives $B_{3}\in\mathcal C(4)$.  The 
equation $S\cdot A_{1}=2B_{1}+2B_{2}+B_{3}$ shows $B_{2}\in\mathcal 
C(4)$.  Finally, we note that $S''\cdot 
A=2B_{1}+6B_{2}+4B_{3}+2B_{4}+2B_{5}$, which shows that $B_{5}\in\mathcal 
C(4)$ as well.  Hence $\overline{\mathcal M}(4,3)_{6}=0$.

\begin{figure}[!h]
			$$\overset{\text{\includegraphics{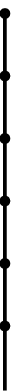}}}{C_{1}}\quad
			\overset{\text{\includegraphics{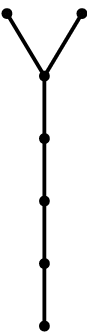}}}{C_{2}}\quad
			\overset{\text{\includegraphics{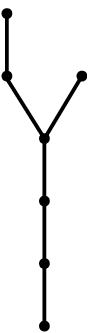}}}{C_{3}}\quad
			\overset{\text{\includegraphics{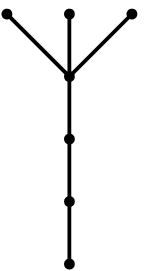}}}{C_{4}}\quad
			\overset{\text{\includegraphics{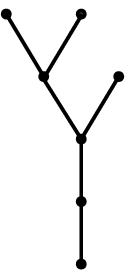}}}{C_{5}}\quad
			\overset{\text{\includegraphics{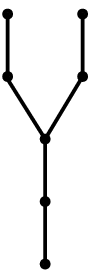}}}{C_{6}}\quad
			\overset{\text{\includegraphics{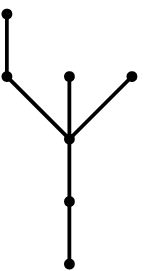}}}{C_{7}}$$
			$$\overset{\text{\includegraphics{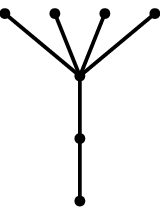}}}{C_{8}}\quad
			\overset{\text{\includegraphics{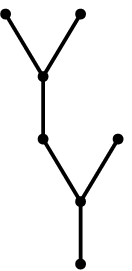}}}{C_{9}}\quad
			\overset{\text{\includegraphics{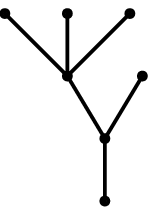}}}{C_{10}}\quad
			\overset{\text{\includegraphics{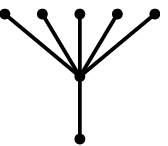}}}{C_{11}}
			$$
			\caption{Trees with 7 vertices}
			\label{Fi:trees7}
			\end{figure}
$\mathcal M_{7}$ is generated over $\mathbb Q$ by the eleven trees 
$C_{1},\ldots,C_{11}$ listed in Figure \ref{Fi:trees7}.  Note that 
$C_{1}$ and $C_{2}$ lie in $\mathcal C(4)$ and that 
$C_{4},C_{7},C_{8},C_{10},C_{11}\in\mathcal V(3)$, which leaves 
$C_{3},C_{5},C_{6},$ and $C_{9}$.  We have 
\begin{align}S\cdot B_{1}=2C_{1}+2C_{2}+2C_{3}&\implies C_{3}\in\mathcal 
C(4)\notag
\\S'\cdot A_{1}=2C_{1}+2C_{3}+C_{6}&\implies C_{6}\in\mathcal C(4)\notag\\
S'''\cdot A=2C_{2}+2C_{5}&\implies C_{5}\in\mathcal C(4)\notag\\
S\cdot B_{2}=C_{2}+2C_{3}+C_{4}+C_{5}+C_{9}&\implies C_{9}\in\mathcal 
N(4,3)\notag
\end{align}
This establishes that $\overline{\mathcal M}(4,3)_{7}=0$.

Lastly we tackle $\mathcal M_{8}$, which is generated over $\mathbb 
Q$ by the unrooted trees $D_{1},\cdots,D_{23}$ given in Figure 
\ref{Fi:trees8}. Apparently $D_{1},D_{2},D_{3},D_{4}\in\mathcal C(4)$ 
and 
$$D_{4}, D_{8}, D_{9}, D_{12}, D_{14} ,D_{15}, D_{16},D_{17},
D_{19},D_{21},D_{22},D_{23}\in\mathcal V(3)\,,$$ leaving us to deal 
with $D_{5},D_{6},D_{7},D_{10},D_{11},D_{13},D_{18},D_{20}$.  Toward 
that end we observe
\begin{align}S\cdot C_{1}=2D_{1}+2D_{2}+2D_{3}+D_{5}&\implies D_{5}\in\mathcal 
C(4)\notag\\
S'\cdot B_{1}=2D_{1}+2D_{3}+2D_{7}&\implies D_{7}\in\mathcal 
C(4)\notag\\
S''''\cdot A=2D_{2}+2D_{10}&\implies D_{10}\in\mathcal C(4)\notag\\
S'''''\cdot A=2D_{3}+2D_{11}&\implies D_{11}\in\mathcal C(4)\notag\\
S'\cdot B_{2}=D_{2}+2D_{5}+D_{8}+D_{10}+D_{13}&\implies 
D_{13}\in\mathcal N(4,3)\notag\\
S'''\cdot A_{1}=2D_{2}+2D_{6}+D_{13}&\implies D_{6}\in\mathcal 
N(4,3)\notag\\
S\cdot C_{2}=D_{2}+2D_{3}+D_{4}+D_{6}+D_{10}+D_{18}&\implies 
D_{18}\in\mathcal N(4,3)\notag\\
S'''\cdot A_{2}=2D_{10}+D_{14}+D_{18}+D_{20}&\implies 
D_{20}\in\mathcal N(4,3)\notag
\end{align}
showing that $\overline{\mathcal M}(4,3)_{8}=0$ and completing the 
proof.

\begin{figure}[!h]
						$$\overset{\text{\includegraphics{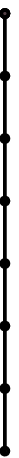}}}{D_{1}}\quad
						\overset{\text{\includegraphics{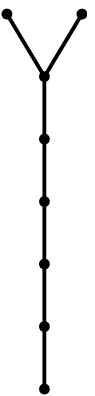}}}{D_{2}}\quad
						\overset{\text{\includegraphics{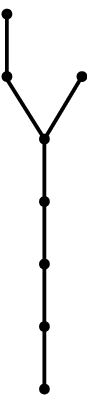}}}{D_{3}}\quad
						\overset{\text{\includegraphics{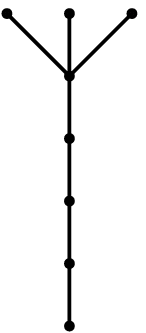}}}{D_{4}}\quad
						\overset{\text{\includegraphics{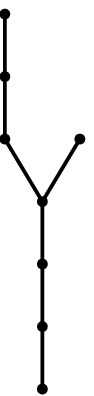}}}{D_{5}}\quad
						\overset{\text{\includegraphics{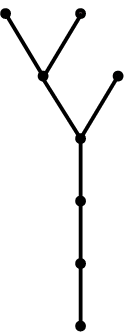}}}{D_{6}}\quad
						\overset{\text{\includegraphics{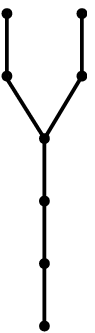}}}{D_{7}}\quad
						\overset{\text{\includegraphics{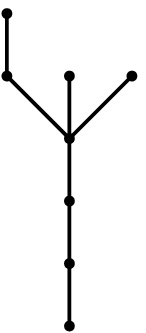}}}{D_{8}}
						$$
						$$
						\overset{\text{\includegraphics{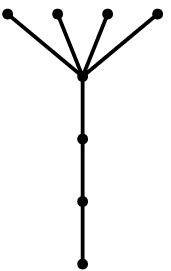}}}{D_{9}}\quad
						\overset{\text{\includegraphics{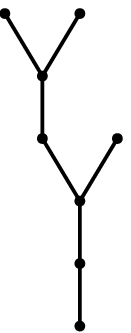}}}{D_{10}}\quad
						\overset{\text{\includegraphics{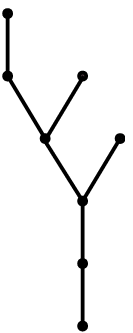}}}{D_{11}}\quad
						\overset{\text{\includegraphics{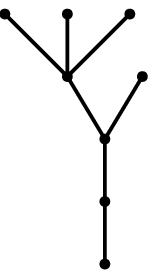}}}{D_{12}}\quad
						\overset{\text{\includegraphics{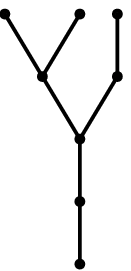}}}{D_{13}}\quad
						\overset{\text{\includegraphics{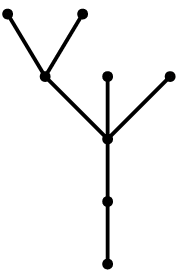}}}{D_{14}}
						$$
						$$
						\overset{\text{\includegraphics{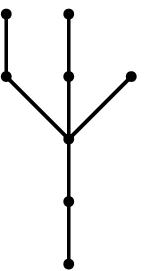}}}{D_{15}}\quad
						\overset{\text{\includegraphics{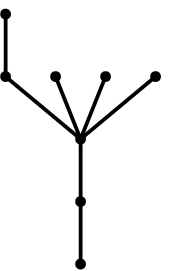}}}{D_{16}}\quad
						\overset{\text{\includegraphics{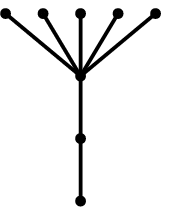}}}{D_{17}}\quad
						\overset{\text{\includegraphics{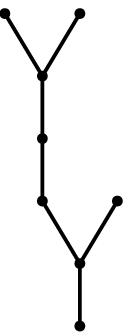}}}{D_{18}}\quad
						\overset{\text{\includegraphics{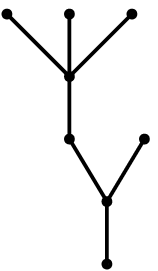}}}{D_{19}}\quad
						\overset{\text{\includegraphics{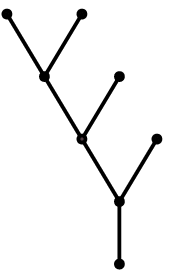}}}{D_{20}}
						$$
						$$
						\overset{\text{\includegraphics{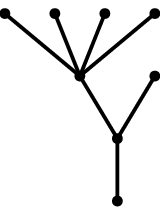}}}{D_{21}}\quad
						\overset{\text{\includegraphics{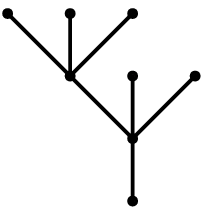}}}{D_{22}}\quad
						\overset{\text{\includegraphics{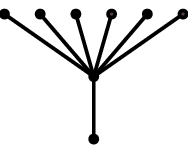}}}{D_{23}}
						$$
						\caption{Trees with 8 vertices}
						\label{Fi:trees8}
						\end{figure}
  \end{proof}
  
\subsection{Questions About the Tree Quotient Modules} The proofs in
sections \ref{cube} and \ref{quad} raise interesting questions about
the tree quotient modules $\overline{\mathcal M}(r,e)$.  For example,
we established in Propositions \ref{3infty} and \ref{43} that
$\overline{\mathcal M}(3,\infty)_{3}=\overline{\mathcal
M}(3,\infty)_{4}=0$.  In fact, the author can prove a far stronger
statement which shows that the tree quotient module
$\overline{\mathcal M}(3,\infty)$ is quite small:

\begin{theo}\label{M3van} $\overline{\mathcal M}(3,\infty)_{m}=0$ 
for $m\ge3$, i.e., $$\overline{\mathcal M}(3,\infty)=\overline{\mathcal 
M}(3,\infty)_{1}\oplus\overline{\mathcal M}(3,\infty)_{2}\,,$$ each 
of these two summands having vector space dimension 1 over $\mathbb Q$.
    \end{theo}

    \noindent The proof will not be given here as it seems to have no 
    implications for the Jacobian Conjecture.

We also established that $\overline{\mathcal M}(4,3)_{m}=0$ for
$5\le m\le8$. One can use the same methods to prove the vanishing 
of $\overline{\mathcal 
M}(4,3)_{m}$ for some larger values of $m$.  So we ask:
    
    \begin{quest}\label{Q2} Is $\overline{\mathcal 
M}(4,3)_{m}=0$ 
    for $m\ge5$?
	\end{quest}
\noindent Of course, an affirmative answer would (seemingly) not
resolve any additional cases of the Symmetric Jacobian Conjecture.
However, an affirmative answer to following question certainly would:

\begin{quest}\label{Q3} Let $r$ be a positive integer. Does there exist a 
positive 
integer $M_{r}$ such that $\overline{\mathcal 
M}(r,4)_{m}=0$ when $M_{r}+1\le m\le2M_{r}$?
\end{quest}

\noindent Or one could ask the weaker question:

\begin{quest}\label{Q3a} Let $r$ be a positive integer. Does there exist a 
positive 
integer $M_{r}$ such that $\overline{\nu}_{m}=0$ in $\overline{\mathcal 
M}(r,4)_{m}$ (see Definition \ref{defnu}) when $M_{r}+1\le m\le2M_{r}$?
\end{quest}
	
\noindent Or one could ask the stronger question:

\begin{quest}\label{Q4} Let $r$ be a positive integer. Does $\overline{\mathcal 
M}(r,4)$ have finite rank as a $\mathbb Q$-vector space?
	\end{quest}
	
	\noindent which is equivalent to asking if $\overline{\mathcal 
	M}(r,4)_{m}=0$ for $m>>0$. It is obvious that the proof of 
	Theorem \ref{d=2,JH4=0Sy} can be mimicked to show that:
	
\begin{theo}\label{solved} Let $r$ be a positive integer.  Assume Question 
\ref{Q3a} (or \ref{Q3}, or \ref{Q4})
has an affirmative answer for $r$, and let $F=X-H$ be a polynomial map
of symmetric homogeneous type with $H$ cubic and $(JH)^{r}=0$.  Then
$F$ is invertible with
$$F^{-1}=X+N^{(1)}+N^{(2)}+N^{(3)}+\cdots+N^{(M_{k})}\,.$$
In particular, the degree of $F^{-1}$ is $\le2M_{k}+1$.
    \end{theo}
    
    Whence, in light of Theorem \ref{SymRed}:
    
\begin{theo} If Question \ref{Q3a} has an affirmative answer for all 
positive integers $r>>0$, then the Jacobian Conjecture is true.
	\end{theo}

	\subsection{The Cubic Symmetric $JH^{4}=0$ Case}\label{CubeQuad}
	
	Questions \ref{Q3} and \ref{Q3a} can be resolved by computer
	algorithm for any fixed $r$, subject to time/space
	limitations.  A computer program has been written and run by
	Li-Yang Tan which appears to resolve the cubic symmetric
	$JH^{4}=0$ case of the Jacobian Conjecture \cite{L-Y}.  The result is
	intriguing.  The program shows $\overline{\mathcal
	M}(4,4)_{m}=0$ for $m=8,9,10,11,12,14$.  Curiously,
	$\overline{\mathcal M}(4,4)_{13}\ne0$ but rather has rank one.
	However the vector $\overline{\nu}_{13}$ (see Definition
	\ref{defnu}) is zero in $\overline{\mathcal M}(4,4)_{13}$.
	Thus we have $\overline{\nu}_{m}=0$ for $m=8,\ldots,14$, so
	the $JH^{4}=0$ case is solved, by Theorem \ref{solved}.  We
	state the theorem thus proved by computer:
	
	\begin{theo} \label{d=3,JH4=0Sy} Let $F=X-H$ be a polynomial map having symmetric
	Jacobian matrix, with $H$ cubic homogeneous and $(JH)^{4}=0$.
	Then $F$ is invertible with
	$$F^{-1}=X+N^{(1)}+N^{(2)}+N^{(3)}+N^{(4)}+N^{(5)}+N^{(6)}+N^{(7)}\,.$$
	In particular, the degree of $F^{-1}$ is $\le15$.
	\end{theo}

  \subsection{Ideal Membership Theorems}
  
  In \cite{W6} the author formulated certain ideal membership
  questions, some of which can be answered in light of the results of
  this paper.  Theorems \ref{JH3=0SymAlt} and \ref{d=2,JH4=0Sym}
  below, which were announced in \cite{W6}, are strengthenings of
  Theorems \ref{JH3=0Sym} and \ref {d=2,JH4=0Sy}, respectively.  In 
  these theorems, the ring $\mathcal  R_{n,[d]}^{\text{\rm sym}}$ is 
  the $\mathbb Q$-algebra generated by the formal coefficients 
  (indeterninates) $c^{q}$ of the formal homogeneous polynomial 
  $P=\sum_{|q|=d+1}c^{q}X^{q}$ of degree $d$.  Here 
  $q=(q_{1},\ldots,q^{n})\in\mathbb N^{n}$, $|q|=q_{1}+\cdots+q_{n}$, 
  and $X_{1}^{q_{1}}=_{n}X^{q_{n}}\cdots X^{q}$.  The reader is 
  referred to \cite{W6} for further explanation of the notation.

  \begin{theo} \label{JH3=0SymAlt} Let $F=X-H$ be the formal degree
  $d\ge2$ polynomial map of symmetric homogeneous type in dimension
  $n$.  In other words $H=\del P$ where $P$ is as above.  Let
  $\mathcal I$ be the ideal in $\mathcal R_{n,[d]}^{\text{\rm sym}}$
  generated by the coefficients of $(JH)^{3}$.  Then all coefficients
  $d^{q}$ of $Q^{(m)}$ for $m\ge3$ (hence all $d^{q}$ with
  $|q|=m(d-1)+2$ with $m\ge3$) are in $\mathcal I$.
  \end{theo}
  
  \begin{proof}  This is immediate from Theorem \ref{JH3=0Sym}, 
  taking $K=\mathcal  R_{n,[d]}^{\text{\rm sym}}/{\mathcal I}$. (This 
  is an advantage to allowing $K$ to be any $\mathbb Q$-algebra.)
      \end{proof}  In similar fashion, the following theorems results 
      from Theorems \ref {d=2,JH4=0Sy} and \ref {d=3,JH4=0Sy}:

\begin{theo} \label{d=2,JH4=0Sym} Let $F=X-H$ be the formal degree $2$
polynomial map of symmetric homogeneous type in dimension $n$, and let
$\mathcal I$ be the ideal in $\mathcal R_{n,[2]}^{\text{\rm sym}}$
generated by the coefficients of $(JH)^{4}$.  Then all coefficients
$d^{q}$ of $Q^{(m)}$ for $m\ge5$ (i.e. $|q|=2m+2$ for $m\ge5$) are in
$\mathcal I$.
\end{theo}

\begin{theo} \label{d=2,JH4=0Sym} Let $F=X-H$ be the formal degree $3$
polynomial map of symmetric homogeneous type in dimension $n$, and let
$\mathcal I$ be the ideal in $\mathcal R_{n,[3]}^{\text{\rm sym}}$
generated by the coefficients of $(JH)^{4}$.  Then all coefficients
$d^{q}$ of $Q^{(m)}$ for $m\ge8$ (i.e. $|q|=m+2$ for $m\ge8$) are in
$\mathcal I$.
\end{theo}\vskip 20pt

\noindent{\bf Acknowledgments: } The author would like to thank his
colleague John Shareshian for his help with numerous computations.  He
also wishes to thank Li-Yang Tan, currently an undergraduate majoring
in Mathematics and Computer Science at Washington University, for
writing the program which yielded the result given in Section
\ref{CubeQuad}.

\renewcommand{\theequation}{\thesection.\arabic{equation}}
\renewcommand{\therema}{\thesection.\arabic{rema}}
\setcounter{equation}{0}
\setcounter{rema}{0}

\bibliographystyle{amsplain}
\bibliography{Refs}

\noindent{\small \sc Department of Mathematics, Washington University in St.
Louis,
St. Louis, MO 63130 } {\em E-mail}: wright@math.wustl.edu

\end{document}